\documentclass[a4paper, 12pt, leqno]{amsart}
\usepackage{mathrsfs}
\usepackage{amsmath}
\usepackage{amscd}
\usepackage{amssymb}
\usepackage{amsthm}

\numberwithin{equation}{section}

\def\Gg{{\mathfrak{g}}}

\def\BC{{\mathbb{C}}}
\def\BF{{\mathbb{F}}}
\def\BN{\mathbb{N}}

\def\BQ{{\mathbb{Q}}}
\def\BZ{{\mathbb{Z}}}
\def\CA{{\mathcal A}}
\def\CB{{\mathcal B}}
\def\CC{{\mathcal C}}
\def\CD{{\mathcal D}}

\def\CN{{\mathcal N}}

\def\lrl{\underset {LR}\leq}
\def\llr{\lrl}

\title[Kazhdan-Lusztig Basis  and A Geometric Filtration
]{Kazhdan-Lusztig Basis  and A Geometric Filtration of  an affine
Hecke Algebra, II}
\author{ Nanhua XI}
\address{Institute of Mathematics\\
Chinese Academy of Sciences\\
Beijing, 100080\\
China } \email{nanhua@math.ac.cn}
\thanks{N. Xi was partially supported by Natural Sciences Foundation
of China (No. 10671193).}

\begin{document}
\begin{abstract}
An affine Hecke algebras can be realized as an equivariant K-group
of the corresponding Steinberg variety. This gives rise naturally to
some two-sided ideals  of the affine Hecke algebra by means of the
closures of nilpotent orbits of the corresponding Lie algebra. In
this paper we will show that  the two-sided ideals are in fact the
two-sided ideals of the affine Hecke algebra defined through
two-sided cells of the corresponding affine Weyl group after the
two-sided ideals are tensored  by $\BQ$. This proves a weak form of
a conjecture of Ginzburg proposed in 1987.
\end{abstract}
\maketitle \setcounter{section}{-1}
\section{Introduction}
\label{sec:Intro} Let $H$ be an affine Hecke algebra over the ring
$\BZ[v,v^{-1}]$ of  Laurent polynomials in an indeterminate $v$ with
integer coefficients. The affine Hecke algebra has a Kazhdan-Lusztig
basis. The basis has many remarkable properties and play an
important role in representation theory. Also, Kazhdan-Lusztig and
Ginzburg give a geometric realization of $H$, which is the key of
Kazhdan-Lusztig's proof for Deligne-Langlands conjecture on
classification of irreducible modules of affine Hecke algebras over
$\BC$ at non-root of 1. This geometric construction of $H$ has  some
two-sided ideals defined naturally by means of the nilpotent variety
of the corresponding Lie algebra. The two-sided ideals form a nice
filtration of the affine Hecke algebra. In [G2] Ginzburg conjectured
that the two-sided ideals are in fact the two-sided ideals of the
affine Hecke algebra defined through two-sided cells of the
corresponding affine Weyl group, see also [L5, T2]. The conjecture
is known to be true for the trivial  nilpotent orbit \{0\} (see
Corollary 8.13 in [L5] and Theorem 7.4 in [X1]) and for type A [TX].
Other evidence is showed in [L5, Corollary 9.13]. We will prove  the
two kinds of two-sided ideals coincide after they are tensored by
$\BQ$ (see Theorem 1.5 in section 1 ). This proves a weak form of
Ginzburg's conjecture.

\section{Affine Hecke algebra}

\noindent{\bf 1.1.} Let $G$ be a simply connected simple algebraic
group over the complex number field $\BC$. The Weyl group $W_0$ acts
naturally on the character group $X$ of a maximal tours of $G$. The
semidirect product $W=W_0\ltimes X$ with respect to the action is
called an (extended) affine Weyl group. Let $H$ be the associated
Hecke algebra over the ring $\CA=\BZ[v,v^{-1}]$ ($v$ an
indeterminate) with parameter $v^2$. Thus $H$ has an $\CA$-basis
$\{T_w\ |\ w\in W\}$ and its multiplication is defined by the
relations $(T_s-v^2)(T_s+1)=0$ if $s$ is a simple reflection and
$T_wT_u=T_{wu}$ if $l(wu)=l(w)+l(u)$, here $l$ is the length
function of $W$.
\def\gb{\mathfrak{b}}

\vskip.3cm \noindent{\bf 1.2.}
Let $\mathfrak{g}$ be the Lie algebra
of $G$, $\mathcal{N}$ the nilpotent cone of $\Gg$ and $\CB$  the
variety of all Borel subalgebras of $\Gg$. The Steinberg variety $Z$
is the subvariety of $\CN\times \CB\times\CB$ consisting of all
triples $(n,\gb,\gb'), \ n\in \gb\cap \gb'\cap\CN, \ \gb,\gb'\in
\CB$. Let $\Lambda=\{(n,\gb)\ | \ n\in\CN\cap \gb, \gb\in\CB\}$ be
the cotangent bundle of $\CB$. Clearly $Z$ can be regarded as a
subvariety of $\Lambda\times\Lambda$ by the imbedding $Z\to
\Lambda\times\Lambda$, $(n,\gb,\gb')\to (n,\gb,n,\gb').$ Define a
$G\times C^*$-action on $\Lambda$ by $(g,z):\ (n,\gb)\to
(z^{-2}\text{ad}(g)n,\text{ad}(g)\gb)$. Let $G\times C^*$ acts on
$\Lambda\times\Lambda$ diagonally, then  $Z$ is a $G\times
C^*$-stable subvariety of $\Lambda\times\Lambda$. For $1\le i<j\le
3$, let $p_{ij}$ be the projection from
$\Lambda\times\Lambda\times\Lambda$ to its $(i,j)$-factor. Note that
the restriction of  $p_{13}$ gives rise to a proper morphism
$p_{12}^{-1}(Z)\cap p_{23}^{-1}(Z)\to Z$. Let
$K^{G\times\mathbb{C}^*}(Z)=K^{G\times\mathbb{C}^*}(\Lambda\times\Lambda;Z)$
be the Grothendieck group of the category of
$G\times\mathbb{C}^*$-equivariant coherent sheaves on
$\Lambda\times\Lambda$ with support in $Z$. We define the
convolution product
$$*:
K^{G\times\mathbb{C}^*}(Z)\times K^{G\times\mathbb{C}^*}(Z)\to
K^{G\times\mathbb{C}^*}(Z),$$
$$\mathscr{F}*\mathscr{G}=(p_{13})_*(p_{12}^*\mathscr{F}\otimes_{\mathcal
{O}_{\Lambda\times\Lambda\times\Lambda}}p_{23}^*\mathscr{G}),$$
where $\mathcal {O}_{\Lambda\times\Lambda\times\Lambda}$ is the
structure sheaf of ${\Lambda\times\Lambda\times\Lambda}$. This
endows with $K^{G\times\mathbb{C}^*}(Z)$ an associative algebra
structure over the representation ring $R_{G\times \BC^*}$ of
$G\times \BC^*$. We shall regard the indeterminate $v$ as the
representation $G\times\BC^*\to \BC^*,\ (g,z)\to z$. Then
$R_{G\times \BC^*}$ is identified with $\CA\otimes_{\BZ} R_G$. In
particular, $K^{G\times\mathbb{C}^*}(Z)$ is an $\CA$-algebra.
Moreover, as $\CA$-algebras, $K^{G\times\mathbb{C}^*}(Z)$ is
isomorphic to the Hecke algebra $H$, see [G1, KL2] or [CG, L5]. We
shall identify $K^{G\times\mathbb{C}^*}(Z)$ with $H$.

\vskip.3cm \noindent{\bf 1.3.} Let $\CC$ and $\CC'$ be two
$G$-orbits in $\CN$. We say that $\CC\le\CC'$ if $\CC$ is in the
closure of $\CC'$. This defines a partial order on the set of
$G$-orbits in $\CN$.  Given a locally closed $G$-stable subvariety
of $\CN$, we set $Z_Y=\{(n,\gb,\gb')\in Z\ |\ n\in Y\}.$

If $Y$ is closed, then the inclusion $i_Y:\ Z_Y\to Z$ induces a map
$(i_Y)_*:\ K^{G\times\mathbb{C}^*}(Z_Y)\to
K^{G\times\mathbb{C}^*}(Z)$ (see [G1, KL2]). The image $H_Y$ of
$(i_Y)_*$ is in fact a two-sided ideal of
$K^{G\times\mathbb{C}^*}(Z)$ (see [L5, Corollary 9.13]), which is
generated by  $G\times \BC^*$-equivariant sheaves supported on
$Z_Y$. It is conjectured that this ideal is spanned by elements in a
Kazhdan-Lusztig basis (see [G2, L5, T2]).

\vskip.3cm \noindent{\bf 1.4.}
 Let $C_w=v^{-l(w)}\sum_{y\le w}P_{y,w}(v^2)T_y$,
where $P_{y,w}$ are the Kazhdan-Lusztig polynomials. Then the
elements $C_w(\ w\in W)$ form an $\CA$-basis of $H$, called a
Kazhdan-Lusztig basis of $H$. Define $w\llr u$ if $a_w\ne 0$ in the
expression $hC_uh'=\sum_{z\in W}a_zC_z\ (a_z\in\CA)$  for some
$h,h'$ in $H$. This defines a preorder on $W$. The corresponding
equivalence classes are called two-sided cells and the preorder
gives rise to a partial order $\llr$ on the set of two-sided cells
of $W$. (See [KL1].) For an element $w$ in $W$ and a two-sided cell
$c$ of $W$ we shall write $w\llr c$ if $w\llr u$ for some
(equivalent to for any) $u$ in $c$.

Lusztig established a bijection between the set of $G$-orbits in
$\CN$ and the set of two-sided cells of $W$ (see [L4, Theorem 4.8]).
Lusztig's bijection preserves the partial orders we have defined,
this was conjectured by Lusztig and verified by Bezrukavnikov (see
[B, Theorem 4 (b)]). Perhaps this bijection is at the heart of the
theory of cells in affine Weyl groups, many deep results are related
to it. Now we can state the main result of this paper.

\vskip.3cm \noindent{\bf Theorem 1.5.} Let $\CC$ be a $G$-orbit in
$\CN$ and $c$ the two-sided cell of $W$ corresponding to $\CC$ under
Lusztig's bijection. Then  the elements $C_w \ (w\llr c)$  form a
$\BQ[v,v^{-1}]$-basis of $H_{\bar\CC}\otimes_\BZ\BQ$, where $\bar
\CC$ denotes the closure of $\CC$ and $H_{\bar\CC}$ is the image of
the map $(i_{\bar{\CC}})_*:\ K^{G\times\mathbb{C}^*}(Z_{\bar\CC})\to
K^{G\times\mathbb{C}^*}(Z)=H$.

\section{Proof of the theorem}

\noindent{\bf 2.1.} Before proving the theorem we need recall some
results about representations of an affine Hecke algebra. Let
$\mathbf H=\BC[v,v^{-1}]\otimes_{\CA} H$ and for any nonzero complex
number $q$ we set $\mathbf{H}_q=\mathbf{H}\otimes_{\BC[v,v^{-1}]}
\BC$, where $\BC$ is regarded as a $\BC[v,v^{-1}]$-algebra by
specializing $v$ to a square root of $q$.¡¡

For any $G$-stable locally closed subvariety $Y$ of $\CN$ we set
$\mathbf{K}^{G\times\BC^*}(Z_Y)=K^{G\times \BC^*}(Z_Y)\otimes\BC$.
If $Y$ is closed, then the inclusion $i_Y:\ Z_Y\to Z$ induces an
injective map $(i_Y)_*:\
\mathbf{K}^{G\times\mathbb{C}^*}(Z_Y)\hookrightarrow
\mathbf{K}^{G\times\mathbb{C}^*}(Z)=\mathbf{H}$. If $Y$ is a closed
subset of $\CN$, we shall identify
$\mathbf{K}^{G\times\mathbb{C}^*}(Z_Y)$ with the image of $(i_Y)_*$,
which is a two-sided ideal of $\mathbf{H}$. See [KL2, 5.3] or [L5,
Corollary 9.13].

Let $s$ be a semisimple element of $G$ and $n$ a nilpotent element
in $\CN$ such that $\text{ad}(s)n=qn$, where $q$ is in $\BC^*$. Let
$\CB_n^s$ be the subvariety of $\CB$ consisting of the Borel
subalgebras containing $n$ and fixed by $s$. Then the component
group  $A(s,n)=C_G(s,n)/C_G(s,n)^o$ of the simultaneous centralizer
in $G$  of $s$ and $n$ acts on the total complex homology group
$H_*(\CB^s_n)$. Let $\rho$ be a representation of $A(s,n)$ appearing
in the space $H_*(\CB^s_n)$. It is known that if $\sum_{w\in
W_0}q^{l(w)}\ne 0$  then the isomorphism classes of irreducible
representations of $\mathbf{H}_q$  is one to one corresponding to
the $G$-conjugacy classes of all the triples $(s,n,\rho),$ where $
s\in G$ is semisimple, $n\in\CN$ satisfying $\text{ad}(s)n=qn$ and
$\rho$ is an irreducible representation of $A(s,n)$ appearing in
$H_*(\CB_n^s)$. See [KL2,  X3].

\def\bhq{\mathbf{H}_q}
\def\bjc{\mathbf{J}_c}

\vskip.3cm \noindent{\bf 2.2.} From now on we assume that $q$ is not
a root of 1.  Let $L_q(s,n,\rho)$ be an irreducible representation
of $\mathbf{H}_q$ corresponding to the triple $(s,n,\rho)$. Kazhdan
and Lusztig constructed a standard module $M(s,n,q,\rho)$ of $\bhq$
such that $L_q(s,n,\rho)$ is the unique simple quotient of
$M(s,n,q,\rho)$ (see [KL2, 5.12 (b) and Theorem 7.12]).
 We shall write $M_q(s,n,\rho)$
 for $M(s,n,q,\rho)$. The following simple fact will be needed.

 \vskip3mm
(a) Let $\CC$ be a $G$-orbit in $\CN$. Then the image
$H_{\bar{\CC}}$ of $(i_{\bar{\CC}})_*$ acts on $M_q(s,n,\rho)$ and
$L_q(s,n,\rho)$ by zero if $n$ is not in $\bar{\CC}$.

{\it Proof.} Clearly $Y=\bar\CC\cup(\overline{G.n}-G.n)$ is closed.
If $n$ is not in $\bar{\CC}$, then the complement in
$X=\bar\CC\cup\overline{G.n}$ of $Y$ is $G.n$. Recall that
$\mathbf{K}^{G\times\mathbb{C}^*}(Z_{Y'})$ is regarded as a
two-sided ideal of $\mathbf{H}$ for any closed subset $Y'$ of $\CN$
(see 2.1). According to [KL2, 5.3 (c), (d) and (e)], the inclusions
$i:\ Y\hookrightarrow X$ and $j:\
 G.n\hookrightarrow X$ induce an exact sequence of
 $\mathbf{H}$-bimodules
 $$0\to \mathbf{K}^{G\times\BC^*}(Z_Y)\to\mathbf{K}^{G\times\BC^*}(Z_X)\to\mathbf{K}^{G\times\BC^*}(Z_{G.n})\to 0.$$
 Using [KL2, 5.3 (e)] we know the inclusion $k:\ \bar\CC\to Y$ induces an injective $\mathbf{H}$-bimodule
 homomorphism $k_*:\ \mathbf{K}^{G\times\BC^*}(Z_{\bar\CC})\to\mathbf{K}^{G\times\BC^*}(Z_Y)$.   Since $M_q(s,n,\rho)$ is a quotient module of
 $\mathbf{K}^{G\times\BC^*}(Z_{G.n})$ (cf. proof of 5.13 in [KL2]),
 the statement (a) then follows from the exact sequence above.

\def\vp{\varphi}
\def\st{\stackrel}
\def\sc{\scriptstyle}
\vskip.3cm \noindent{\bf 2.3.}
Let $J_c$ be the based ring of a
two-sided cell $c$ of $W$, which has a $\BZ$-basis  \{$t_w\, |\, w
\in c$\}. Let $D_c$ be the set of distinguished involutions in $c$.
For $x,y\in W$, we write $C_xC_y=\sum_{z\in W}h_{x,y,z}C_z,\
h_{x,y,z}\in\CA$. The map
$$\vp_c(C_w)=\displaystyle\sum_{\st{\sc d\in D_c}{\st{\sc u\in W}
{a(d)=a(u)}}}h_{w,d,u}t_u,\qquad w\in W,$$ defines an $\CA$-algebra
homomorphism $H\to J_c\otimes_{\BZ}\CA$, where $a:W\to \BN$ is the
$a$-function defined in [L1]. The homomorphism $\vp_c$ induces a
$\BC$-algebra homomorphism $\vp_{c,q}:\ \bhq\to\bjc=J_c\otimes_{\BZ}
\BC$. If $E$ is a $\bjc$-module, through $\vp_{c,q}$, $E$ gets an
$\bhq$-module structure, which will be denoted by $E_q$. See [L2,
L3].

Let $\CC$ be the nilpotent orbit corresponding to $c$. According to
[L4, Theorems 4.2 and  4.8], the map $E\to E_q$ defines a bijection
between the isomorphism classes of simple $\bjc$-modules and the
isomorphism classes of standard modules  $M_q(s,n,\rho)$ with $n$ in
$\CC$.

Now we start to  prove Theorem 1.5.

\vskip.3cm \noindent{\bf 2.4.} We first show that  $H_{\bar{\CC}}$
is contained in the two-sided ideal $H^{\le c}$ spanned by all $C_w
 \ (w\llr c)$.

Let $\CC=G.n$ and recall that $H_{\bar\CC}$ stands for the image of
$(i_{\bar\CC})_*: \ K^{G\times\BC^*}(Z_{\bar\CC})\to
K^{G\times\BC^*}(Z)=H$. If $H_{\bar\CC}$ was not contained in the
$\CA$-submodule $H^{\le c}$ of $H$, we could find $x\in W$ such that
$x\not\llr c$ and $C_x$ appears in $H_{\bar\CC}$. (We say that $C_x$
appears in $H_{\bar\CC}$ if there exists an element $\sum_{w\in
W}a_wC_w$ ($a_w\in \CA$) in $H_{\bar{\CC}}$ such that $a_x\ne 0$.)
Choose $x\in W$ such that $C_x$ appears in $H_{\bar\CC}$ and $x$ is
highest with respect to the preorder $\llr$ and to $H_{\bar\CC}$ in
the following sense: whenever $C_w$ appears in $H_{\bar\CC}$, then
either $w$ and $x$ are in the same two-sided cell or $x\not\llr w$.
Let $c'$ be the two-sided cell containing $x$. We then have
$c'\not\llr c$.

Choose an element $h=\sum_{w\in W}a_wC_w$ ($a_w\in \CA$)  in
$H_{\bar{\CC}}$ such that $h_{c'}=\sum_{w\in c'}a_wC_w$ is nonzero.
Then we have $\varphi_{c'}(h)=\varphi_{c'}(h_{c'})$.

We claim that $\varphi_{c'}(h_{c'})$ is nonzero. Let $u\in c'$ be
such that $a_u$ has the highest degree (as a Laurent polynomial in
$v$) among all $a_w, \ w\in c'$. Let $d$ be the distinguished
involution such that $d$ and $u$ are in the same left cell. Then we
know that for any distinguished involution $d'$, the degree
$h_{w,d',u}$  is less than
 the degree of $h_{u,d,u}$ if either $w\ne u$ or $d'\ne d$ (see
[L2, Theorems 1.8 and 1.10]). Thus  the degree of  $a_wh_{w,d',u}$
is less than the degree of $a_uh_{u,d,u}$ if either $w\ne u$ or
$d'\ne d$. Hence $\varphi_{c'}(h_c')$ is nonzero.

 Clearly,
there are only finitely many $q$ such that $\varphi_{c',q}(h_{c'})$
is zero after specializing $v$ to a square root of $q$. According to
Theorem 4 in [BO], the  ring $\mathbf{J}_{c'}$ is semisimple, that
is, its Jacobson radical is zero. So we can find a nonzero $q$ in
$\BC$ with infinite order and a simple $\mathbf{J}_{c'}$ module $E'$
such that $\varphi_{c',q}(h)=\varphi_{c',q}(h_{c'})$ is nonzero and
acts on $E'$ is not by zero.

According to [L4, Theorems 4.2 and 4.8], $E'_q$ is isomorphic to
certain standard module $M_q(s',n',\rho)$ with $n'$ in the nilpotent
orbit $\CC'$ corresponding to $c'$. Since $c'\not\llr c$, $\CC'$ is
not in the closure of $\CC$ (see [B, Theorem 4 (b)]), by 2.2 (a),
the image $H_{\bar\CC}$ of $(i_{\bar{\CC}})_*$
 acts on $E'_q$ by zero. This contradicts that $\varphi_{c',q}(h)$ acts on
$E'$ is not by 0. Therefore $H_{\bar{\CC}}$ is contained in the
two-sided ideal $H^{\le c}$.

\vskip.3cm \noindent{\bf 2.5.} In this subsection all tensor
products are over $\BZ$ except other specifications are given.

 Now
we show that $H^{\le c}\otimes \BQ$ is equal to
$H_{\bar{\CC}}\otimes \BQ$. If $\CC$ is regular, then $\bar\CC$ is
the whole nilpotent cone and the corresponding two-sided cell $c$
contains the neutral element $e$, in this case, both $H_{\bar{\CC}}$
and $H^{\le c}$ are the whole Hecke algebra.

We use induction on the partial order $\llr$ in the set of all
two-sided cells of $W$. Assume that for all $c'$ with $c\llr c'$ and
$c'\ne c$, we have $H_{\bar{\CC'}}\otimes \BQ=H^{\le c'}\otimes
\BQ$, where $\CC'$ is the nilpotent orbit corresponding to $c'$.

We need to show $H_{\bar{\CC}}\otimes \BQ=H^{\le c}\otimes \BQ$. Let
$c'$ be a two-sided cell different from $c$ such that $c\llr c'$ but
there is no two-sided cell $c''$  between $c$ and $c'$, i.e. no
$c''$ such that $c\llr c''\llr  c'$ and  $c\ne c''\ne c'$.

Let $\BF$ be an algebraic closure of $\BC(v)$. We first show that
$\BF\otimes_\CA H_{\bar \CC}=\BF\otimes_\CA H^{\le c}$. Assume this
was not true. Note that $\BF$ is isomorphic to $\BC$
(non-canonically), so we can apply the results in [KL2].  By 2.4 and
induction hypothesis, there would exist $w\in c$ such that $C_w$ is
contained in $\BF\otimes_\CA H_{\bar\CC'}$ but not in
$\BF\otimes_\CA H_{\bar \CC}$. By the choice of $\CC'$, we know that
$\BF\otimes_\CA H_{\bar\CC'-\CC'}$ is the sum of $\BF\otimes_\CA
H_{\bar\CC}$ and $\BF\otimes_\CA H_{\bar\CC''}$ for some nilpotent
orbits $\CC''$ with $\bar\CC\not\subseteq\bar\CC''$ and
$\bar\CC''\not\subseteq\bar\CC$ (see [KL2, 5.3(e)]). Since
$\bar\CC\not\subseteq\bar\CC''$, by  [B, Theorem 4 (b)] we know that
$C_w$ is not in $H^{\le c''}$, where $c''$ is the two-sided cell
corresponding to $\CC''$. By 2.4 we see that $\BF\otimes_{\CA}
H_{\bar\CC''}$ does not contain $C_w$. Thus the image in
$M_{\CC'}=\BF\otimes_\CA K^{G\times \BC^*}(Z_{\CC'})= \BF\otimes_\CA
H_{\bar \CC'}/\BF\otimes_\CA H_{\bar \CC'-\CC'}$ of $C_w$ is
nonzero. According to [KL2, Corollary 5.9], each nonzero element in
$\BF\otimes_\CA H_{\bar\CC'} \setminus \BF\otimes_\CA H_{\bar
\CC'-\CC'}$ acts on $M_{\CC'}$ by nonzero. The argument for [KL2,
Proposition 5.13] implies that each nonzero element in $M_{\CC'}$
would have nonzero image in some standard quotient module
 of $M_{\CC'}$. Thus $C_w$ acts on some standard quotient module $M_{v^2}(s,n',\rho)$ of $M_{\CC'}$
 by nonzero, where $n'\in\CC'$.

 The homomorphism $\vp_{c'}$ (see 2.3 for definition) induces a homomorphism $\varphi_{c',v^2}:\
\BF\otimes_\CA H\to \BF\otimes J_{c'}$. According to [L4, Theorems
4.2 and 4.8], $M_{v^2}(s,n',\rho)$ is isomorphic
 to certain $E_{v^2}$ (defined similar to the $E_q$ in 2.3) for some simple $\BF\otimes J_{c'}$-module $E$.
 Note that $h_{w,d,u}\ne 0$ implies that $u\llr w$. Since $\CC'$ is not in the closure of $\CC$
 and $w$ is in the two-sided
$c$ corresponding to $\CC$,
 using  [L4, Theorem
4.8] and [B, Theorem 4 (b)], we see that $\varphi_{c',v^2}(C_w)=0$.
Then $C_w$ acts $M_{v^2}(s,n',\rho)$ by zero. This leads to a
contradiction. So we have $\BF\otimes_\CA H_{\bar
\CC}=\BF\otimes_\CA H^{\le c}$.

Thus for each $w\in c$, we can find a nonzero $a\in\BF$ such that
$aC_w$ is in $H_{\bar\CC}$. Clearly, we must have $a\in\CA$.  Now we
show that
 $\mathbf{K}^{G\times \BC^*}(Z_Y)$ is a free
$\BC[v,v^{-1}]$-module for any $G$-stable locally closed subvariety
$Y$ of $\CN$. According to [KL2, 5.3] we may assume that $Y$ is a
nilpotent orbit $\CC$. It is enough to show that the completion of
$\mathbf{K}^{G\times \BC^*}(Z_\CC)$ at any semisimple class in
$G\times C^*$ is free over $\BC[v,v^{-1}]$. Using [KL2, 5.6] it is
enough to show that the right hand side of 5.6(a) in [KL2] is free.
This follows from [KL2, (l3)], the assumption there is satisfied by
[KL2, 4.1]. Using [KL2, 5.3] we know that as free
$\BC[v,v^{-1}]$-modules $H_{\bar\CC'}\otimes\BC$ is a direct sum of
$H_{\bar\CC}\otimes\BC$ and $\mathbf{K}^{G\times
\BC^*}(Z_{\bar\CC'-\bar\CC})$. By assumption,  $H_{\bar\CC'}\otimes
\BQ=H^{\le c'}\otimes \BQ$, thus $H_{\bar\CC'}\otimes \BQ$ is a free
$\BQ[v,v^{-1}]$-module and contains $C_w$. These imply that if
$aC_w\in H_{\bar\CC}$ for some nonzero $a\in\CA$ then $C_w\in
H_{\bar\CC}\otimes \BC$. Therefore we can find a nonzero complex
number $a$ such that $aC_w$ is in $H_{\bar\CC}$. Obviously we have
$a\in\BZ$. Thus $H^{\le c}\otimes \BQ$ is contained in
$H_{\bar{\CC}}\otimes \BQ$. By 2.4 we then have  $H^{\le c}\otimes
\BQ=H_{\bar{\CC}}\otimes \BQ$. Theorem 1.5 is proved.

\section{Some comments}

\noindent{\bf 3.1.} If one can show that $K^{G\times\BC^*}(Z_\CC)$
is a free $\BZ$-module for any nilpotent orbit  $\CC$, then the
argument in 2.5 shows that the image of $(i_{\bar{\CC}})_*$ in
$H=K^{G\times\mathbb{C}^*}(Z)$ contains $H^{\le c}$, where $c$ is
the two-sided cell corresponding to $\CC$. Thus Ginzburg's
conjecture would be proved. In fact, it seems that one can expect
more. More precisely, it is likely the following result is true.

\vskip.3cm

 (a) $K^{G\times\BC^*}(Z_\CC)$ is a free $\CA$-module and
$K_1^{G\times\BC^*}(Z_\CC)$=0 for all nilpotent orbit $\CC$. (We
refer to [CG, section 5.2] and [Q] for the definition of the functor
$K_i^G.$) \vskip.3cm

 If (a) is true, then we also have

(b) The map  $(i_{\bar{\CC}})_*: K^{G\times\BC^*}(Z_{\bar\CC})\to
K^{G\times\BC^*}(Z)$ is injective.

\vskip.3cm

\def\vp{\varphi}
We explain some evidences for (a) and prove it for $G=GL_n(\BC),$ $
Sp_4(\BC)$ and type $G_2$. Let $N$ be a nilpotent element in $\CC$
and $\CB_N$ be the variety of Borel subalgebras of $\Gg$ containing
$N$. By the Jacobson-Morozov theorem, there exists a homomorphism
$\varphi:\ SL_2(\BC)\to G$ such that $d\varphi\begin{pmatrix}
0&1\\0&0\end{pmatrix}=N$. For $z$ in $\BC^*$, let
$d_z=\begin{pmatrix} z&0\\0&z^{-1}\end{pmatrix}$. Following Kazhdan
and Lusztig [KL2, 2.4], we define $Q_N=\{(g,z)\in G\times\BC^*\ |\
\text{ad}(g)N=z^2N\}.$  Then $Q_N$ is a closed. Let $x=(g,z)\in Q_N$
act on $(G\times \BC^*)\times \CB_N\times\CB_N$ by
$x(y,\gb,\gb')=(yx^{-1}, \text{ad}(g)\gb,\text{ad}(g)\gb')$. Then
$Z_\CC$ is isomorphic to the quotient space $Q_N\backslash((G\times
\BC^*)\times \CB_N\times\CB_N)$. Thus we have
$K_i^{G\times\BC^*}(Z_\CC)=K_i^{Q_N}(\CB_N\times\CB_N)$ (see [KL2,
5.5] and [Th1, Prop. 6.2]). It is known that  $Q_\varphi=\{(g,z)\in
G\times C^*\ |\ g\vp(x) g^{-1}=\vp(d_zxd_z^{-1})\text{ for all }x\in
SL_2(\BC)\}$ is a maximal reductive subgroup of $Q_N$ (see [KL2, 2.4
(d)]). So we have
$K_i^{Q_N}(\CB_N\times\CB_N)=K_i^{Q_\vp}(\CB_N\times\CB_N)$ (see
[CG, 5.2.18]).

Let $P$ be the parabolic subgroup of $G$ associated to $N$ (see
[DLP, 1.12]). Then we know that the intersection $\CB_{N,\mathscr
O}$ of $\CB_N$ with any $P$-orbit $\mathscr O$ on $\CB$ is smooth.
The torus $\CD=\{\vp(d_z)\ |\ z\in\BC^*\}$ is a subgroup of $P$ and
acts on $\CB_{N,\mathscr O}$,  and $\CB_{N,\mathscr O}$ is a vector
bundle over the $\CD$-fixed point set $\CB_{N,\mathscr O}^{\CD}$
(see [DLP, 3.4 (d)]). Since the action of $Q_\vp$ on
$\CB_{N,\mathscr O}$ commutes with the action of $\CD$,  according
to [BB],  this vector bundle is isomorphic to a $Q_\vp$-stable
subbundle of $T(\CB_{N,\mathscr O})|_{\CB_{N,\mathscr O}^{\CD}}$,
where $T(\CB_{N,\mathscr O})$ is the  tangent bundle of
$\CB_{N,\mathscr O}$. Thus the vector bundle is $Q_\vp$-equivariant,
so that the computation of $K_i^{Q_\vp}(\CB_N\times\CB_N)$ is
reduced to the computation of $K_i^{Q_\vp}(\CB_{N,\mathscr
O}^\CD\times\CB_{N,\mathscr O'}^\CD)$ for various $P$-orbits
$\mathscr{O,\ O'}$ on $\CB$ (see Theorems 2.7 and  4.1 in [Th1], or
Theorems 5.4.17 and 5.2.14 in [CG]). Note that
$C_\vp=\{g\vp(d_z^{-1})\,|\, (g,z)\in Q_\vp\}$ is a maximal
reductive subgroup of the centralizer $C_G(N)$ of $N$ (see [BV,
2.4]) and the map $(g,z)\to (g\vp(d_z^{-1}),z)$ define an
isomorphism from $Q_\vp$ to $C_\vp\times \BC^*$. Thus we have
$K_i^{Q_\vp}(\CB_{N,\mathscr O}^\CD\times\CB_{N,\mathscr
O'}^\CD)=K_i^{C_\vp\times \BC^*}(\CB_{N,\mathscr
O}^\CD\times\CB_{N,\mathscr O'}^\CD)$. Now the factor $\BC^*$ and
the group $\CD$ act on $\CB_{N,\mathscr O}^\CD\times\CB_{N,\mathscr
O'}^\CD$ trivially, we therefore have $K_i^{Q_\vp}(\CB_{N,\mathscr
O}^\CD\times\CB_{N,\mathscr O'}^\CD)=K_i^{C_\vp}(\CB_{N,\mathscr
O}^\CD\times\CB_{N,\mathscr O'}^\CD)\otimes R_{\BC*}$ (see [CG,
(5.2.4)], the argument there works for higher $K$-groups). Note that
we have identified $R_{\BC*}$ with $\CA=\BZ[v,v^{-1}]$. Thus the
statement (a) is equivalent to the following one.

\vskip.3cm

(c) $K_i^{C_\vp}(\CB_{N,\mathscr O}^\CD\times\CB_{N,\mathscr
O'}^\CD)$ is a free $\BZ$-module for $i=0$ and is 0 for $i=1$.

\vskip.3cm

The statement (c) seems much  easier to access. The variety
$\CB_{N,\mathscr O}^\CD$ and its fixed point set $\CB_{N,\mathscr
O}^{s,\CD}$ by any semisimple element $s$ in $C_\vp$ are smooth and
have good homology properties. See[DLP].

\vskip.3cm

\noindent{\bf Proposition 3.2.} The statement (a) is true for
$GL_n(\BC),\ Sp_4(\BC)$ and type $G_2$. In particular, Ginzburg's
conjecture is true in these cases.

{\it Proof.} We only need to prove statement (c). For $G=GL_n(\BC)$,
we know that $\CB_{N,\mathscr O}^\CD$ has an $\alpha$-partition into
subsets which are affine space bundles over the flag variety $\CB'$
of $C_\vp$ (see Theorem 2.2 and 2.4 (a) in [X2]). In this case, (a)
is true since we are reduced to compute
$K_i^{C_\vp}(\CB'\times\CB')$ (cf. [CG, Lemma 5.5.1] and the
argument for [L6, Lemma 1.6]). For $G=Sp_4(\BC)$ or type $G_2$, we
know that $\CB_{N,\mathscr O}^\CD$ is either empty or the flag
variety of $C_\vp$ if $N$ is not subregular (see Prop. 4.2 (i) and
section 4.4 in [X2]. In this case, we are also reduced to compute
$K_i^{C_\vp}(\CB'\times\CB')$ ({\it loc.cit}), so (a) is true. If
$N$ is subregular, then $\CB_N$ is a Dynkin curve and it is easy to
see that $\CB_{N,\mathscr O}^\CD$ is either a projective line or a
finite set (see Prop. 4.2 (ii) and section 4.4 in [X2] for a
computable description of $\CB_N$). The computation for
$K_i^{C_\vp}(\CB_{N,\mathscr O}^\CD\times\CB_{N,\mathscr O'}^\CD)$
is easy, it is a free $\BZ$-module for $i=0$, see 4.3 (b) and 4.4 in
[X2], and is 0 for $i=1$ (since this is true for projective line and
a finite set). Then we have the same conclusion for
$K_i^{C_\vp}(\CB_{N}\times\CB_{N})\ ({\it loc.cit})$. The
proposition is proved.

\vskip3mm Remark:  For $GL_n(\BC)$, this proposition   also provide
an another proof for the main result of [TX], where results in [T1]
are used.

\vskip3mm
 \noindent{\bf Proposition 3.3.}  Assume that $C_\vp$  is
connected. Then

\noindent (a) $K^{C_\vp}(\CB_N\times\CB_N)$ is a free $\BZ$-module.

\noindent (b) $K^{Q_\vp}(\CB_N\times\CB_N)$ is a free $\CA$-module.
That is, $K^{G\times \BC*}(Z_{{G.N}})$ is a free $\CA$-module.

{\it Proof.} Let $T$ be a maximal torus. According to [Th2, (1.11)],
we have split monomorphism $K^{C_\vp}(\CB_N\times\CB_N)\to
K^{T}(\CB_N\times\CB_N)$. Similar to the argument for [L6, Lemma
1.13 (d)], we see that $K^{T}(\CB_N\times\CB_N)$ is a free
$R_T$-module. (a) follows.

The reasoning for (b) is similar since $Q_\vp$ is isomorphic to
$C_\vp\times C^*$ and the monomorphism
$K^{Q_\vp}(\CB_N\times\CB_N)\to K^{T\times\BC^*}(\CB_N\times\CB_N)$
is split. The proposition is proved.

\vskip3mm

Remark: If $G=GL_n(\BC)$, then all $C_\vp$ are connected and have
simply connected derived group. In this case
$K^{Q_\vp}(\CB_N\times\CB_N)$ is a free $R_{Q_\vp}$-module since
$R_{Q_\vp}=R_{C_\vp}\otimes \CA$ and  $R_{T\times \BC^*}$ is a free
$R_{C_\vp}\otimes \CA$-module. Combining this, subsection 2.4 and
the argument in subsection 2.5 we obtain a different proof for the
main result in [TX].

\vskip3mm

\noindent{\bf 3.4.} The $K$-groups $K^{F}(\CB_N)$ and
$K^{F}(\CB_N\times\CB_N)$ are important in representation theory of
affine Hecke algebras for $F$ being $Q_\vp,\ C_\vp$ or a torus of
$Q_\vp$ (see [KL2, L6]). For the nilpotent element $N$, in [L4,
10.5] Lusztig conjectures there exists a finite $C_\vp$-set $Y$
which plays a key role in understanding the based ring of the
two-sided cell corresponding to $G.N$. It seems that as
$R_{C_\vp}$-modules $K^{C_\vp}(Y)$ and $K^{C_\vp}(Y\times Y)$ are
isomorphic to $K^{C_\vp}(\CB_N)$ and $K^{C_\vp}(\CB_N\times \CB_N)$
respectively. Let  $X=\CB_N\text{ or }\CB_N\times\CB_N$. In view of
[L4, 10.5] one may hope to find a canonical $\BZ$-basis of
$K^{C_\vp}(X)$ and a canonical $\CA$-basis of $K^{Q_\vp}(X)$ in the
spirit of [L5, L6]. Moreover, there should exist a natural bijection
between the  elements of the canonical basis of $K^{F}(\CB_N\times
\CB_N)$ ($F=C_\vp$ or $Q_\vp$) and the elements of the two-sided
cell corresponding to $G.N$.

\vskip3mm
 {\bf Acknowledgement:} I thank Professor G. Lusztig for very
helpful correspondences and for telling me the argument for the
freeness of $\mathbf{K}^{G\times C^*}(Z_\CC)$ over $\BC[v,v^{-1}]$.
I am grateful to Professor T. Tanisaki for a helpful correspondence
and to Professor Jianzhong Pan for a helpful conversation.

\bibliographystyle{unsrt}

\begin{thebibliography}{}
\bibitem[B]{B}
Bezrukavnikov, R.: Perverse sheaves on affine flags and nilpotent
cone of the Langlands dual group, math.RT/ 0201256v4.
\bibitem[BB]{BB}
Bialynicky-Birula, A.: Some theorems on actions of algebraic groups,
Ann. Math. \textbf{98} (1973), 480-497.
\bibitem[BO]{BO}
 Bezrukavnikov, R.,  Ostrik, V.: On tensor categories
attached to cells in affine Weyl groups, II, in "Representation
theory of algebraic groups and quantum groups", pp.101--119, Adv.
Stud. Pure Math., 40, Math. Soc. Japan, Tokyo, 2004.
\bibitem[BV]{BV}
Barbash, D., Vogan, D.: Unipotent representation of complex
semisimple groups, Ann. Math. \textbf{121} (1985), 41-110.
\bibitem[CG]{CG}
Chriss, N., Ginzburg, V.: Representation theory and complex
geometry, Birkh\"auser Boston, Inc., Boston, MA, 1997.
\bibitem[DLP]{DLP}
DeConcini, C., Lusztig, G., Procesi, C.: Homology of the zero-set of
a nilpotent vector field on a flag manifold, Jour. Amer. Math. Soc.
\textbf{1} (1988), 15-34.
\bibitem[G1]{G1}
Ginzburg, V.: Lagrangian construction of representations of Hecke
algebras, Adv. in Math. \textbf{63} (1987), 100--112.
\bibitem[G2]{G2}
Ginzburg, V.: Geometrical aspects of representation theory,
Proceedings of the International Congress of Mathematicians, Vol.
\textbf{1,} (Berkeley, Calif., 1986), 840--848, Amer. Math. Soc.,
Providence, RI, 1987.
\bibitem[KL1]{KL1}
Kazhdan, D., Lusztig, G.: Representations of Coxeter groups and
Hecke algebras,  Invent. Math. \textbf{53}  (1979), 165--184.
\bibitem[KL2]{KL2}
Kazhdan, D., Lusztig, G.: Proof of the Deligne-Langlands conjecture
for Hecke algebras, Invent. Math. \textbf{87}  (1987), 153--215.
\bibitem[L1]{L:cell1}  Lusztig, G.:  Cells in affine Weyl groups, in
``Algebraic groups and related topics", Advanced Studies in Pure
Math., vol. \textbf{6}, Kinokuniya and North Holland, 1985, pp.
255-287.
\bibitem[L2]{L:cell2}  Lusztig, G.:  Cells in affine Weyl groups, II,
J. Alg. \textbf{109} (1987), 536-548.
\bibitem[L3]{L:cell3}  Lusztig,
G.:  Cells in affine Weyl groups, III,  J. Fac. Sci. Univ. Tokyo
Sect. IA Math. \textbf{34}  (1987), 223--243.
\bibitem[L4]{L:cell4}
Lusztig, G.: Cells in affine Weyl groups. IV, J. Fac. Sci. Univ.
Tokyo Sect. IA Math. \textbf{36}  (1989), 297--328.
\bibitem[L5]{L5}
Lusztig, G.: Bases in equivariant $K$-theory, Represent. Theory
\textbf{2}  (1998), 298--369 (electronic).
\bibitem[L6]{L6} Lusztig,
G.: Bases in equivariant $K$-theory, II, Represent. Theory
\textbf{3} (1999), 281--353 (electronic).
\bibitem[Q]{Q}
Quillen, D.: Higher algebraic K-theory I, in {\it Higher
$K$-Theory,} Lecture Notes in Math. \textbf{341} (1972), 85-147.
\bibitem[T1]{T1}
Tanisaki, T.: Hodge modules, equivariant $K$-theory and Hecke
algebras, Publ. Res. Inst. Math. Sci. \textbf{23}  (1987), 841--879.
\bibitem[T2]{T2}
Tanisaki, T.: Representations of semisimple Lie groups and
$D$-modules, Sugaku expositions \textbf{4}  (1991), 43--61.
\bibitem[Th1]{Th1}
Thomason, R.W.: Algebraic K-theory of group scheme actions, in {\it
Algebraic Topology and Algebraic K-Theory,} Ann. of Math. Studies
\textbf{113} (1987), 539--563.
\bibitem[Th2]{Th2} Thomason, R.W.:
Equivariant algebraic vs. topological K-homology Atiyah-Segal-Style,
Duke Math. J. \textbf{56} (1988), 589--636.
\bibitem[TX]{TX}
Tanisaki, T. and Xi, N.: Kazhdan-Lusztig Basis  and A Geometric
Filtration of an affine Hecke Algebra, Nagoya Math. J. {\bf 182}
(2006), 285-311.
\bibitem[X1]{Xi1}
Xi, N.: Representations of affine Hecke algebras, Lecture Notes in
Mathematics, \textbf{1587}. Springer-Verlag, Berlin, 1994.
\bibitem[X2]{Xi2}
Xi, N.: A partition of the Springer fibers $\mathcal B_N$ for type
$A_{n-1}, B_2, G_2$ and some applications, Indag. Mathem., N.S.,
\textbf{10}(2) (1999), 307-320.
\bibitem[X3]{Xi3}
Xi, N.: Representations of affine Hecke algebras and based ring of
affine Weyl groups, J. Amer. Math. Soc. {\bf 20} (2007)£¬211-217.
\end{thebibliography}

\end{document}